\documentclass[letterpaper, 11 pt, draftcls, onecolumn]{ieeeconf}

\IEEEoverridecommandlockouts                              

\usepackage{amsmath} 
\usepackage{amssymb}  
\usepackage{color}
\usepackage{graphicx}
\usepackage{booktabs} 

\newcommand{\A}{\mathcal{A}}
\newcommand{\U}{\mathcal{U}}
\newcommand{\C}{\mathcal{C}}
\newcommand{\E}{\mathcal{E}}

\newcommand{\DAM}{[\partial \A]_{-}}
\newcommand{\Int}{\mathsf{int}}
\newcommand{\ee}{\epsilon}
\newcommand{\dee}{\,\mathrm{d}}

\newcommand{\dist}{\mathsf{dist}}

\newtheorem{remark}{Remark}
\newtheorem{definition}{Definition}
\newtheorem{theorem}{Theorem}
\newtheorem{example}{Example}
\newtheorem{proposition}{Proposition}
\newtheorem{corollary}{Corollary}

\title{Recursive feasibility of continuous-time model predictive control without stabilising constraints$^{\star}$} 

\author{Willem Esterhuizen$^{1}$, Karl Worthmann$^{2}$ and Stefan Streif$^{1}$
\thanks{$^{1}$Willem Esterhuizen \& Stefan Streif, Technische Universit\"{a}t Chemnitz, Automatic Control and System Dynamics Laboratory, Chemnitz, Germany;
        {\tt\scriptsize [willem.esterhuizen, stefan.streif]@etit.tu-chemnitz.de}
        }%
\thanks{$^{2}$Karl Worthmann, Institute for Mathematics, Technische Universit\"{a}t Ilmenau, Ilmenau, Germany;
        {\tt\scriptsize karl.worthmann@tu-ilmenau.de}}%
\thanks{$^{\star}$ This work has been accepted in IEEE Control Systems Letters, 10.1109/LCSYS.2020.3001514. ©2020 IEEE. Personal use of this material is permitted. Permission from IEEE must be obtained for all other uses, in any current or future media, including reprinting/republishing this material for advertising or promotional purposes, creating new collective works, for resale or redistribution to servers or lists, or reuse of any copyrighted component of this work in other works.}%
}

\begin{document}

\maketitle
\thispagestyle{empty}
\pagestyle{empty}

\begin{abstract}
	We consider sampled-data Model Predictive Control (MPC) of nonlinear continuous-time control systems. %
	We derive sufficient conditions to guarantee recursive feasibility and asymptotic stability without stabilising costs and/or constraints. %
	Moreover, we present formulas to explicitly estimate the required length of the prediction horizon based on the concept of (local) cost controllability. %
	For the linear-quadratic case, %
	cost controllability can be inferred from standard assumptions. %
	In addition, we extend results on the relationship between the horizon length and the distance of the initial state to the boundary of the viability kernel %
	from the discrete-time to the continuous-time setting.
\end{abstract}


\section{Introduction}\label{sec:introduction}



\noindent Two key aspects in Model Predictive Control (MPC) are asymptotic stability and recursive feasibility, see, e.g.\ \cite{RawlMayn18}. %
The latter property refers to existence of a solution to the optimisation problem recursively invoked in the MPC algorithm. %
Both properties are often ensured by constructing suitable terminal costs and constraints, which are then added to the optimisation problem. 
However, recently, researchers have thoroughly analysed MPC schemes without stabilising constraints, i.e.\ so-called unconstrained MPC, w.r.t.\ their stability behaviour, see, e.g. \cite{NeviPrim97,JadbYu01,GrimMess05,GrunPann10} for references dealing with discrete-time systems and \cite{REBLE2012,Worthmann2014} for extensions to continuous-time systems and the relation to the former. See \cite{cannon_2015} for a study of how these two approaches are linked. The main motivation behind this conceptual shift is simplicity of the approach and its reduced numerical complexity, which explains its pre-dominant use in industry, cp.\ the discussion in \cite[Section~7.4]{GrunPann17}. %

However, a rigorous treatment of recursive feasibility has only been done for discrete-time systems, see~\cite{GrunPann17,BOCCIA2014,boccia_2014_linear}, while in the continuous-time case only particular cases were studied, see \cite{FaulBonv15b,FaulBonv15}. 
On the one hand, only qualitative statements (sufficiently large prediction horizon) were given, but without any estimates on the required length. On the other hand, more restrictive assumptions were made, i.e.\ a combination of local controllability (linearized system is controllable) and reachability, see~\cite[assumption~(ii) of Proposition~1]{FaulBonv15}. Indeed, we implicitly show the second assumption en passant in our line of reasoning. %
In the references \cite{GrunPann17,BOCCIA2014,boccia_2014_linear}, the authors show that sublevel sets of the finite-horizon value function can be made recursively feasible assuming so-called cost-controllability~\cite{Coron_2019}, i.e., the same property as typically invoked for showing asymptotic stability. To this end, a sufficiently long prediction horizon is required, which is 
quantified in dependence of the problem data. %
In this paper, we derive a continuous-time analogon of these results: 
We establish a condition that allows one to determine the length of the prediction horizon such that recursive feasibility and asymptotic stability are 
guaranteed. To this end, a thorough investigation of the interplay between open- and closed-loop control is of key importance. %
In addition, we consider the Linear-Quadratic (LQ) case in detail, showing that the viability kernel~\cite{DeDona_siam,aubin2009viability} is closely related to the sublevel sets of the finite-horizon value function. 

The outline of the paper is as follows: In Section~\ref{sec:setting}, we present the setting, i.e.\ nonlinear continuous-time systems with state and input constraints, and the MPC algorithm. Section~\ref{sec:rec-stab} presents our main result, i.e.\ a sufficient condition for recursive feasibility and asymptotic stability. Then, we focus on the LQ case in Section~\ref{sec:LQR} 
before our findings 
are illustrated by an example in Section~\ref{sec:example} and conclusions are drawn in Section~\ref{sec:conclusion}.

\textbf{Notation}: $\mathbb{R}_{\geq 0}$ and $\mathbb{R}_{>0}$ to the set of nonnegative and positive real numbers resp. %
For $x \in \mathbb{R}^n$, $|x|$ denotes its Euclidean norm %
while $\Vert F \Vert := \max_{|x| = 1} |Fx|$ denotes the induced matrix norm of $F \in \mathbb{R}^{m \times n}$. %
For a set $S \subseteq \mathbb{R}^q$ with $q \geq n$, $\partial S$ and $\mathsf{proj}_{\mathbb{R}^n}(S)$ refer to its boundary and projection onto $\mathbb{R}^n$ resp. A continuous function $\eta: \mathbb{R}_{\geq 0} \rightarrow \mathbb{R}_{\geq 0}$ belongs to class $\mathcal{K}_{\infty}$ if it is strictly increasing, unbounded, and $\eta(0) = 0$. The space $L_{\mathrm{loc}}^{\infty}(\mathbb{R}_{\geq 0},\mathbb{R}^m)$ denotes the set of Lebesgue-measurable functions $u: \mathbb{R}_{\geq 0} \rightarrow \mathbb{R}^m$ that are locally essentially bounded. Given two sets $S_1\subset\mathbb{R}^n$ and $S_2\subset\mathbb{R}^n$, $\mathrm{dist}(S_1,S_2):=\inf_{x\in S_1} \mathrm{dist}(x,S_2)$ where $\mathrm{dist}(x,S)$ is the Euclidean distance from $x\in\mathbb{R}^n$ to $S\subset \mathbb{R}^n$.

\section{Problem setting}\label{sec:setting}

\noindent We consider the continuous-time nonlinear system:
\begin{align}
\dot{x}(t) & = f(x(t),u(t))	 \label{eq:system-1}
\end{align}
with initial condition $x(0) = x_0$, where $x(t)\in\mathbb{R}^n$ and $u(t)\in\mathbb{R}^m$ are the state and the control at time 
$t \in \mathbb{R}_{\geq 0}$ resp. Let the map $f:\mathbb{R}^n\times\mathbb{R}^m\rightarrow\mathbb{R}^n$ be continuous and locally Lipschitz continuous with respect to its first argument~$x$ on $\mathbb{R}^n \setminus \{0\}$. Then, for a given control function $u \in L_{\mathrm{loc}}^{\infty}(\mathbb{R}_{\geq 0},\mathbb{R}^m)$, there exists a unique solution of the initial value problem, which is denoted by $x(t) = x(t;x_0,u)$, $t \in I_{x_0,u}$, where $I_{x_0,u}$ denotes its maximal interval of existence, see for example \cite[App. C]{sontag2013}. 
{Using the set $\mathcal{E}\subseteq\mathbb{R}^n\times\mathbb{R}^m$, we impose the constraints} 
\begin{align}
	(x(t),u(t)) & \in \mathcal{E} \qquad\forall\, t\in[0,\infty). \label{eq:system-3}
\end{align}
Moreover, the sets $U(x) := \{u\in\mathbb{R}^m : (x,u)\in\mathcal{E}\}$, $x \in X$, and $X:=\mathsf{proj}_{\mathbb{R}^n}(\mathcal{E}) = \{x: U(x) \neq \emptyset 
\}$ are introduced. With $x_0\in X$ and $T\in\mathbb{R}_{> 0}$, we let $\mathcal{U}_T(x_0)$ (resp.\ $\mathcal{U}_{\infty}(x_0)$) denote the set of all $u\in L_{\mathrm{loc}}^{\infty}(\mathbb{R}_{\geq 0},\mathbb{R}^m)$ such that the solution exists and satisfies the constraint~\eqref{eq:system-3} on $[0,T]$ (resp.\ for all $t \geq 0$). We will refer to the set $\mathcal{U}_T(x_0)$ (resp. $\mathcal{U}_{\infty}(x_0)$) as the set of \emph{admissible control functions} w.r.t.\ the initial {value~$x_0$} and the horizon~$T$ (resp.\ the infinite horizon).

Next, we define the {\emph{viability kernel}~\cite{aubin2009viability} (also called admissible set in, e.g.,\ \cite{DeDona_siam})} as the set of states, for which an admissible control exists on $[0,\infty)$.
\begin{definition}[Viability Kernel]\ 
\label{def:admissible-set}
	The {viability kernel} 
	is defined by $\A := \{x_0 \in \mathbb{R}^n : \mathcal{U}_{\infty}(x_0) \neq \emptyset \}$.
\end{definition}


Let $ \bar{x} \in X$ be a controlled equilibrium, that is, there exists a $\bar{u} \in U( \bar{x})$ such that $f(\bar{x},\bar{u}) = 0$ holds. Our goal is to use MPC to asymptotically steer the state to $\bar{x}$ while satisfying the state and control constraints~\eqref{eq:system-3}. To that end, we define a finite-horizon cost functional
\[
	J_T(x_0,u):= \int_0^T \ell(x(s;x_0,u), u(s) )\, \dee s,
\]
with continuous \textit{stage cost} $\ell:\mathbb{R}^n \times \mathbb{R}^m \rightarrow \mathbb{R}_{\geq 0}$, and consider the finite-horizon {Optimal Control Problem (OCP)}
\begin{equation}\label{eq:finit-horizon-OCP}
	\inf_{u\in\mathcal{U}_T(x_0)}J_T(x_0,u).
\end{equation}
The value function $V_T: X\rightarrow \mathbb{R}_{\geq 0} \cup \{+\infty\}$ is defined as $V_T(x_0) := \inf_{u\in\mathcal{U}_T(x_0)} J_T(x_0,u)$. By convention, we will say that $V_T(x_0) = \infty$ when $\mathcal{U}_T(x_0) = \emptyset$. Moreover, we use the abbreviation $V_T^{-1}[0,C] := \{ x \in X | V_T(x) \leq C\}$.

We impose the following assumption on $\ell$, which is, {e.g., satisfied for quadratic stage cost $\ell(x,u) := x^T Q x + u^T R u$ with positive definite weighting $Q \in \mathbb{R}^{n \times n}$ if $\bar{x} = 0$, $\bar{u} = 0$}.
\begin{description}
	\item[(A1)] 
	{Let there be {$\mathcal K_{\infty}$-functions $\underline{\eta}, \overline{\eta}$} satisfying 
		\begin{align*}
			\underline{\eta}(|x -  \bar{x}|) \leq \ell^{\star}(x) \leq \overline{\eta}(|x -  \bar{x}|) \qquad \forall\, x\in X
		\end{align*}
		using the abbreviation $\ell^{\star}(x) := \inf_{u \in U(x)} \ell(x,u)$.}
\end{description}

\noindent Given a time shift $\delta \in\mathbb{R}_{>0}$, a number $N\in\mathbb{N}$ that specifies the length of the prediction (optimization) horizon, and an initial value~$x_0 \in X$, the MPC algorithm is as follows:
\begin{enumerate}
	\item[1.] Set the prediction horizon $T \leftarrow N\delta$, let $p \leftarrow 0$ 
	\item[2.] Measure the current state $\hat{x} = x(p\delta;x_0,u_{\operatorname{MPC}})$
	\item[3.] Find a minimiser $u^{\star} \in \arg\inf_{u\in\mathcal{U}_T(\hat{x})} J_T(\hat{x},u)$
	\item[4.] Implement $u_{\operatorname{MPC}}(t) = u^{\star}(t)$, $t\in[p\delta,(p+1)\delta)$
	\item[5.] Set $p \leftarrow p + 1$ and go to step~2
\end{enumerate}
We assume that a minimiser of problem~\eqref{eq:finit-horizon-OCP} in step~3 exists if $\mathcal{U}_T(x_0) \neq \emptyset$ holds 
to facilitate 
the upcoming analysis, see for e.g. \cite[p.56]{GrunPann17} for a detailed discussion of that assumption.

The MPC algorithm implicitly defines a sampled-data feedback law $\mu_{T,\delta}: [0,\delta) \times X \rightarrow \mathbb{R}^m$, $\mu_{T,\delta}(t,\hat{x}) = u^\star(t)$. Here, $u^\star \in \mathcal{U}_T(\hat{x})$ denotes the solution of the finite-horizon optimal control problem depending on the initial state~$\hat{x}$. This (open-loop) control is applied on the time interval $[0,\delta)$. The resulting \textit{closed-loop} solution at time~$t$ emanating from $x_0$ is denoted by $x_{\mu_{T,\delta}}(t;x_0)$. For practical purposes, the set of admissible control functions is often limited to sampled-data with zero-order hold, e.g., $u^\star(t) \equiv \hat{u} \in \mathbb{R}^m$, $t \in [0,\delta)$. %
Algorithmically, ensuring admissibility of a control function, i.e.\ $u \in \mathcal{U}_T(\hat{x})$, is a non-trivial task. 
Hence, the derivation of sufficient conditions such that only finitely many inequality constraints have to be checked to ensure admissibility is of particular interest, see, e.g.\ \cite[Lemma~1]{DaruWort18} for an example tailored to the non-holonomic robot. 
One of our goals for future research is to alleviate this burden by using the characterization of the viability kernel's boundary based on the so-called theory of barriers~\cite{DeDona_siam}. There it is shown that parts of the boundary are made up of integral curves of the system that satisfy a minimum-like principle, which yields conditions on the control.

The following property is essential in order to ensure well-posedness of the MPC algorithm.
\begin{definition}\label{def:recurs-feas-set}
	A set $\mathcal{S}\subseteq X$ is said to be \emph{forward invariant} w.r.t. the sampled feedback law $\mu_{T,\delta}$ if and only if, for each $\hat{x} \in \mathcal{S}$, the conditions
	\begin{itemize}
		\item $(x(t;\hat{x},u^\star),u^\star(t)) \in \E$ for all $t \in [0,\delta]$ and
		\item $x(\delta;\hat{x},u^\star) \in\mathcal{S}$
	\end{itemize}		
	hold with $u^\star(t) := \mu_{T,\delta}(t;\hat{x})$ depending on $\hat{x}$ and $T$.
\end{definition}

If the initial state is located in a forward-invariant set, then the state and input remain feasible on $[0, \delta]$ under the MPC feedback, and the state is contained in this set at time $\delta$. Then, the OCP to be solved in the third step of the proposed MPC algorithm is feasible in each iteration as can be shown by induction. In conclusion, one gets recursive feasibility of the MPC closed loop. %
Note that a forward-invariant set w.r.t.\ some feedback is, of course, a subset of the viability kernel but, in general, not vice versa. 
However, the viability kernel may contain states that cannot be steered to the origin as shown in the following example.
\begin{example}\label{ExampleSetO}
	Consider the scalar system $\dot{x}(t) = x(t) + u(t)$ with the constraints $|x(t)| \leq 2$ and $|u(t)| \leq 1$. Then, $\partial \A$ is given by the set $\{-1,1\}$. For each $x_0 \in \partial \A$, $u(t) \equiv - x_0$ with state trajectory $x(t;x_0,u) \equiv x_0$ is the only element of $\mathcal{U}_\infty(x_0)$. 
	Hence, we get $V_\infty(x_0) = \infty$ for all $x_0 \in \partial \A$ if the stage cost satisfies Assumption~A1.
\end{example}	

Hence, our goal is to show both forward invariance (and, thus, recursive feasiblity of the MPC closed loop) and asymptotic stability for a (hopefully) large set~$\mathcal{S} \subseteq \mathcal{A}$.


%

\section{Recursive feasibility and asymptotic stability}\label{sec:rec-stab}

To address stability and recursive feasibility of the MPC closed loop, we impose \emph{cost controllability}~\cite{Coron_2019}. 
\begin{description}
	\item[(A2)] There exists $\gamma\in\mathbb{R}_{>0}$ and a neighbourhood~$\mathcal N$ of $\bar{x}$ such that $V_{\infty}(x)\leq \gamma \ell^{\star}(x)$ holds for all $x\in\mathcal{N}\cap X$.
\end{description}
Note that we employ the local version following~\cite{BOCCIA2014}, which assumes cost controllability only in a neighbourhood of the controlled equilibrium~$\bar{x}$. It is possible to further weaken Assumption~A2 by using a growth bound~$\gamma$ depending on the length of the prediction horizon, see, e.g.\ \cite{Wort11} and the references therein for a detailed discussion. Assumption~A2, in fact, holds on arbitrary sublevel sets of the finite-horizon value function. To show this, note that the constant
\begin{equation}\label{eq:lower-bound-on-l}
	M := \inf_{x\in X\setminus\mathcal{N}} \ell^{\star}(x)> 0 
\end{equation}
is well-defined in view of Assumption~A1.
\begin{proposition}
	Let Assumptions~A1 and~A2 hold. Then, for each $C\in\mathbb{R}_{>0}$ such that $\mathcal{N} \cap X \subseteq V_T^{-1}[0,C]$, the following inequality holds with $\beta := \max\{C/M,\gamma\}$
	\begin{equation}\label{ineq:value-level-sets}
		V_T(x) \leq \beta \ell^{\star}(x) \qquad \forall\, x\in V_T^{-1}[0,C]. 
	\end{equation}
\end{proposition}
\begin{proof}
	Using Inequality~\eqref{eq:lower-bound-on-l}, we have $V_T(x) \leq \frac{C}{M} \ell^{\star}(x)$ for all $x\in V_T^{-1}[0,C]\setminus \mathcal{N}$. Moreover, from Assumption A2 we have $V_T(x)\leq \gamma \ell^{\star}(x)$ for all $x\in \mathcal{N}\cap X$. Combining these two inequalities yields the assertion.
\end{proof}
{We impose the following assumption, which can be considered as a \textit{consistency condition} ensuring that low-order terms of the value function can be estimated from below by using the stage cost, cp.~the framework on stage cost design presented in~\cite{Coron_2019} for further insight, see also Section~\ref{sec:LQR}, in which we verify Assumption~A3 in the LQ case.}
\begin{description}
	\item[(A3)] For given prediction horizon~$T > 0$ and $C > 0$, there exists a constant $\bar{C}\in\mathbb{R}_{>0}$ such that {the following inequality holds for all $\delta \in (0,T]$}
	\[
		\delta \ell^{\star}(\hat{x}) \leq \bar{C} V_\delta(\hat{x}) \qquad\forall\; \hat{x} \in V_T^{-1}[0,C].
	\]
\end{description}

We now present the main contribution of the paper: the extension of \cite[Theorem~4]{BOCCIA2014} from discrete to continuous time. Theorem~\ref{thm:main} shows that the (controlled) equilibrium~$\bar{x}$ is asymptotically stable w.r.t.\ the sampled-data feedback law~$\mu_{T,\delta}$ with domain of attraction containing the set $V_T^{-1}[0,C]$. The latter implies recursive feasibility of the MPC closed loop for arbitrary initial values $x_0 \in V_T^{-1}[0,C]$.
\begin{theorem}\label{thm:main}
	Consider the system~\eqref{eq:system-1} 
	and constraint~\eqref{eq:system-3}. %
	For {given $C > 0$, let} Assumptions~A1, A2 with~$\gamma$, and A3 with~$\bar{C}$ hold. %
	Moreover, let $M$ be defined {by}~\eqref{eq:lower-bound-on-l} and $\beta:=\max\{\frac{C}{M},\gamma\}$. %
	Then, for $\delta \in (0,\beta)$, and prediction horizon {$T = N \delta$, $N\in\mathbb{N}$,} satisfying the condition
	\begin{equation}
		\max \Bigg\{ \frac {C}{M\delta}, %
		\bar{C} \left(\frac{\beta}{\delta}\right)^2 \Bigg\} \cdot \left(\frac{\beta}{\beta+\delta}\right)^{N-1} < 1, \label{ineq-2}
	\end{equation}
	the following relaxed Lyapunov inequality holds for all $\hat{x} \in V_T^{-1}[0,C]$ with $\alpha = \alpha_{N,\delta} := \bar{C} (\beta/\delta)^2 (\beta/(\beta + \delta)^{N-1})$:
	\begin{align}\label{ineq:RelaxedLyapunov}
		V_T(x(\delta;\hat{x},\mu_{T,\delta})) \leq V_T(\hat{x}) - (1-\alpha) \int_{0}^{\delta} \ell(s)\,\mathrm{d}s,
	\end{align}
	where $\ell(s) := \ell(x(s;\hat{x},u^\star),u^\star(s))$.
\end{theorem}
\begin{proof} 
	Consider $\hat{x} \in V_T^{-1}[0,C]$, and let $u^{\star}\in\mathcal{U}_T(\hat{x})$ be a solution of the OCP~\eqref{eq:finit-horizon-OCP}, %
	which implies $J_T(\hat{x},u^{\star}) \leq J_T(\hat{x},u)$ for all $u \in \mathcal{U}_T(\hat{x})$. %
	Define $\tilde{x} := x(\delta;\hat{x},u^{\star})$, and let $\tilde{u}^\star$ be the respective solution of~\eqref{eq:finit-horizon-OCP}, %
	i.e.\ $J_T(\tilde{x},\tilde{u}^\star) = V_T(\tilde{x})$.\footnote{Note that this equality also holds ($V_T(\tilde{x}) = \infty$) if there does not exist an admissible solution of the OCP~\eqref{eq:finit-horizon-OCP}.} For any $t \in [0,T]$, we have
	\begin{align}
		V_T(\hat{x}) & = \int_0^{\delta} \hspace*{-0.075cm} \ell(s)\,\mathrm{d}s + \hspace*{-0.025cm} \int_{\delta}^{t} %
		\hspace*{-0.075cm} \ell(s)\,\mathrm{d}s + \hspace*{-0.025cm} \int_t^{T} \hspace*{-0.075cm} \ell(s)\,\mathrm{d}s \label{eq:value-x0}
	\end{align}
	and $V_T(\tilde{x}) = \int_{0}^T \tilde{\ell}(s)\,\mathrm{d}s \leq \int_{\delta}^{t} \ell(s)\,\mathrm{d} s + V_{T-t+\delta}(x(t;\hat{x},u^\star))$ %
	with $\tilde{\ell}(s) := \ell(x(s;\tilde{x},\tilde{u}^\star),\tilde{u}^\star(s))$. Hence, using~\eqref{eq:value-x0} to replace $\int_\delta^t \ell(s)\,\mathrm{d}s$ in the last inequality yields
	\begin{align}\label{eq:vlaue-at-next-step}
		V_T(\tilde{x}) & \leq V_T(\hat{x}) - \hspace*{-0.025cm} \int_0^{\delta} \hspace*{-0.075cm} \ell(s)\,\mathrm{d}s - \hspace*{-0.025cm} %
		\int_t^{T} \hspace*{-0.075cm} \ell(s)\,\mathrm{d}s + V_{T-t+\delta}(x(t))
	\end{align}
	with $x(t) := x(t;\hat{x},u^\star)$. 
	
	Before we proceed, we make the following preliminary considerations: %
	Since the assumption $V_T(\hat{x}) \leq C$ implies $V_{T - t}(x(t;\hat{x},u^\star)) \leq C$ for all $t \in [0,T]$, %
	using Inequality~\eqref{ineq:value-level-sets} yields $V_{T-t}(x(t;\hat{x},u^\star)) \leq \beta \ell^\star(x(t;\hat{x},u^\star)) \leq \beta \ell(t)$. %
	Therefore, we have $V_T(\hat{x}) \leq \int_0^{t} \ell(s)\,\mathrm{d}s + \beta \ell(t)$ and, thus, %
	$\int_{t}^{T} \ell(s)\,\mathrm{d}s \leq \beta \ell(t)$ for all $t \in [0,T]$. %
	For any $\bar{t}\in[0, T - \delta]$, this implies the inequality
	\[
	\int_{\bar{t}}^{\bar{t} + \delta} \left( \int_{t}^T \ell(s)\dee s \right) \mathrm{d} t \leq \beta \int_{\bar{t}}^{\bar{t} + \delta} \ell(t)\,\mathrm{d}t,
	\]
	which is, using $\int_t^T \ell(s)\dee s = \int_{t}^{\bar{t} + \delta} \ell(s) \dee s + \int_{\bar{t} + \delta}^T \ell(s)\dee s $, equivalent to:
	\[
	\int_{\bar{t}}^{\bar{t} + \delta} \hspace*{-0.05cm} \left(\int_{t}^{\bar{t} + \delta} \hspace*{-0.05cm} \ell(s)\dee s \right) \hspace*{-0.05cm}\mathrm{d}t + \delta \int_{\bar{t} + \delta}^T \hspace*{-0.05cm} \ell(s)\dee s \leq \beta \int_{\bar{t}}^{\bar{t} + \delta} \hspace*{-0.075cm} \ell(t) \dee t.
	\]
	Since the first term on the left-hand side is nonnegative, we get $\delta \int_{\bar{t} + \delta}^T \ell(s)\dee s \leq \beta \int_{\bar{t}}^{\bar{t} + \delta} \ell(s) \dee s$ for all $\bar{t} \in[0,T-\delta]$. Thus, setting $\bar{t} = p \delta$ for $p \in \{0,1,\dots,N-1\}$, leads to
	\begin{align*}
	\frac{\delta}{\beta}\int_{(p+1)\delta}^{N\delta} \ell(s) \dee s & \leq \int_{p\delta}^{(p+1)\delta} \ell(s) \dee s.
	\end{align*}
	Adding the term $\int_{(p+1)\delta}^{N\delta} \ell(s) \dee s$ on both sides of this inequality reads as
	\begin{align*}
	\frac{\beta + \delta}{\beta}\int_{(p+1)\delta}^{N\delta} \ell(s) \dee s & \leq \int_{p\delta}^{N\delta}  \ell(s) \dee s.
	\end{align*}
	Iteratively applying this inequality until $p = N-2$ to further estimate the left hand side yields
	\begin{align*}
		\int_{p\delta}^{N\delta}  \ell(s) \dee s 
	& \geq \left(\frac{\beta + \delta}{\beta}\right)^{N - p -1}\int_{(N-1)\delta}^{N\delta} \ell(s) \dee s
	\end{align*}
	for $p \in \{0,1,\ldots,N-1\}$. Then, first taking $p = 0$ and invoking $\hat{x} \in V_T^{-1}[0,C]$ Assumption~A2, and the definition of~$\beta$ afterwards, yields
	\begin{equation}
	\min\{\beta \ell^{\star}(\hat{x}), C\} \geq V_T(\hat{x})\geq \left(\frac{\beta + \delta}{\beta}\right)^{N-1}\int_{T - \delta}^{T} \ell(s) \dee s.\label{eq:bounds-on-V}
	\end{equation}
	Using (the first term of) Condition~\eqref{ineq-2}, we get $\int_{T - \delta}^T \hspace*{-0.075cm} \ell(s)\,\mathrm{d}s < M\delta$, which implies $\int_{T-\delta}^T \ell^{\star}(x(s;\hat{x},u^{\star}))\,\mathrm{d}s < M \delta$. Since
	\begin{equation}\label{ineq:l-star-at-t-hat}
		\delta \left( \inf_{t \in [T-\delta,T]} \ell^{\star}(x(t;\hat{x},u^{\star})) \right) \leq %
		\int_{T-\delta}^T \hspace{0.075cm} \ell^{\star}(x(s;\hat{x},u^{\star}))\,\mathrm{d}s
	\end{equation}
	holds, there exists a time instant $\hat{t} \in [T-\delta,T]$ such that the inequality $\ell^{\star}(x(\hat{t};\hat{x},u^{\star})) < M$ holds. %
	As a consequence, the definition of~$M$ ensures $x(\hat{t};\hat{x},u^{\star})\in\mathcal{N}$, cp.~\eqref{eq:lower-bound-on-l}. %
	Therefore, we have $V_{\infty}(x(\hat{t};\hat{x},u^{\star})) \leq \beta \ell^{\star}(x(\hat{t};\hat{x},u^{\star}))$, which implies\footnote{We stress that the following inequality implies finiteness of $V_T(\tilde{x})$ and, thus, the existence of an admissible control function $\tilde{u} \in \mathcal{U}_T(\tilde{x})$.}
	\[
	V_{T - \hat{t} + \delta}(x(\hat{t};\hat{x},u^{\star})) \leq \beta \ell^{\star}(x(\hat{t};\hat{x},u^{\star})).
	\]
	Note that the definition of~$\hat{t}$ and the line of reasoning used to derive Inequality~\eqref{ineq:l-star-at-t-hat} imply
	\begin{align*}
		\delta \ell^{\star}(x(\hat{t};\hat{x},u^{\star})) \leq \int_{T-\delta}^T \hspace*{-0.075cm} \ell^{\star}(x(s;\hat{x},u^{\star})\,\mathrm{d}s \leq \int_{T-\delta}^T \hspace*{-0.075cm} \ell(s)\,\mathrm{d}s.
	\end{align*}
	Thus, using~\eqref{eq:vlaue-at-next-step} with $t = \bar{t}$ and applying the last two inequalities, we get %
	\begin{align*}
		V_T(\tilde{x}) &\leq V_T(\hat{x}) - \hspace*{-0.0375cm} \int_0^{\delta} \hspace*{-0.15cm} \ell(s)\,\mathrm{d}s - \hspace*{-0.0375cm} \int_{\hat{t}}^{T} \hspace*{-0.15cm} \ell(s)\,\mathrm{d}s + \frac{\beta}{\delta} \int_{T-\delta}^T \hspace*{-0.1cm} \ell(s)\,\mathrm{d}s.
	\end{align*}
	Then, dropping the term $\int_{\hat{t}}^T \ell(s)\,\mathrm{d}s$, $\hat{t}\in[T-\delta,T]$, and using Inequality~\eqref{eq:bounds-on-V}, we get
	\begin{align*}
		V_T(\tilde{x}) \leq V_T(\hat{x}) - \hspace*{-0.025cm} \int_0^{\delta} \hspace*{-0.075cm} \ell(s)\,\mathrm{d}s + \frac{\beta^2}{\delta} \left( \frac{\beta}{\beta + \delta}\right)^{\hspace*{-0.075cm} N-1} \ell^\star(\hat{x}).
	\end{align*}	
	Invoking Assumption~A3, we get the desired relaxed Lyapunov inequality~\eqref{ineq:RelaxedLyapunov} where $\alpha_{N,\delta}$ is defined as in (the second term of) Condition~\eqref{ineq-2}, i.e.\ the assertion.
\end{proof}

Theorem~\ref{thm:main} relates the easily checkable Condition~\eqref{ineq-2} with the Lyapunov-like decrease~\eqref{ineq:RelaxedLyapunov} of the finite-horizon value function, which allows one to conclude asymptotic stability of the equilibrium and can be, thus, considered as a sufficient stability condition, cp.\ \cite{GrunPann10}. %

For the following statement, we require the set of exceptional points defined by $\mathcal{O}:=\lim_{n\rightarrow \infty} V_{\infty}^{-1}[n,\infty)$, cp.\ \cite{BOCCIA2014}. The set~$\mathcal{O}$ consists of all elements of the state space at which the value function may blow up. Hence, entering the set~$\mathcal{O}$ should be avoided. Using this definition, we can derive the following corollary of Theorem~\ref{thm:main}, which is the continuous-time analogon of \cite[Thm. 6]{BOCCIA2014}, which involves compact subsets of the infinite-horizon value function's sublevel sets.
\begin{corollary}\label{cor:compact-subsets}
	Let {Assumptions A1-A3} hold and let $K \subset V_{\infty}^{-1}[0,\infty) \setminus \mathcal{O}$ be a compact set. %
	Then, for a chosen time shift $\delta > 0$ there exists {a} horizon length $\bar{T} = \bar{T}(\delta,K) \in (\delta,\infty)$ such that %
	the origin is asymptotically stable w.r.t. the MPC closed loop with domain of attraction containing the set~$K$.
\end{corollary}
\begin{proof}
	By the same arguments as in \cite{BOCCIA2014}, for any compact set $K\subset V_{\infty}^{-1}[0,\infty)\setminus\mathcal{O}$, %
	there exists a $C<\infty$ such that $K\subseteq V_{\infty}^{-1}[0,C]$. %
	Therefore, from Condition~\eqref{ineq-2} and with a chosen time shift $\delta$, for any initial condition $x_0\in K$ %
	the MPC closed loop is stable and recursively feasible if $N \geq \bar{N}$, i.e.\ $\bar{T} := \delta \bar{N}$, %
	where $\bar{N} = \bar{N}(\delta, K) \in \mathbb{N}$ satisfies
	\begin{equation}
		\bar{N} > \frac{ \max\{ \ln(\delta) + \ln(M) - \ln(C), 2\ln(\delta/\beta) - \ln(\bar{C}) \} }{\ln(\beta) - \ln(\beta + \delta)} + 1.\label{eq:ineq_2}
	\end{equation}
\end{proof}

We stress that Corollary~\ref{cor:compact-subsets} implies recursive feasibility for all initial values~$x_0$ contained in the compact set $K$. We clarify the the importance of considering compact sets $K\subset V_{\infty}^{-1}[0,\infty)\setminus \mathcal{O}$ in the next section.

In the following section, we show that some of the constants used in Theorem~\ref{thm:main} can be determined in the LQ case using the spectrum of the weighting matrices in the stage cost and the solution to the algebraic Riccati equation. In the nonlinear settting, the system at hand has to be considered in detail, see, e.g.\ \cite{WortMehr15}, but is, in general, a non-trivial task.

\section{Constrained Linear-Quadratic Case}\label{sec:LQR}

\noindent In Section~\ref{sec:rec-stab} we presented conditions {rendering 
sublevel sets of the finite-horizon value function forward invariant w.r.t.\ $\mu_{T,\delta}$ 
and ensuring asymptotic stability of the equilibrium.} %
We now consider the LQ case and establish a relationship between the viability kernel, see Definition~\ref{def:admissible-set}, and the level sets of the value function. Then, we specialise Corollary~\ref{cor:compact-subsets} to the LQ case; an important result because there exist algorithms capable of computing compact inner-approximations of viability kernels, see \cite{saintpierre_1994,monnet_2016}. 
Most of the results derived for discrete-time systems in~\cite{BOCCIA2014} carry over to our continuous-time setting in an analogous way. Hence, we only present \textit{novel} aspects in the main text and refer the reader to the appendix for details. 

We focus on linear 
systems, i.e.\ system~\eqref{eq:system-1} given by
\begin{equation}
	f(x(t),u(t)) := Ax(t) + Bu(t) \label{sys:linear}
\end{equation}
with matrices $A\in\mathbb{R}^{n\times n}$ and $B\in\mathbb{R}^{n \times m}$. Moreover, we impose pure control constraints, i.e.\ $u(t) \in U$, $U \subset \mathbb{R}^m$, for all $t \geq 0$ and the following assumptions:
\begin{description}
	\item[(A4)] The pair $(A,B)$ is stabilisable
	\item[(A5)] {The sets $U$, $\E := \{ (x,u) \in \mathbb{R}^{n+m} : g(x,u) \leq 0\}$} %
		with $g \in \mathcal{C}^2(\mathbb{R}^n \times \mathbb{R}^m ) \rightarrow \mathbb{R}^p$ %
		are convex, compact, and contain the origin in their interior. %
		The mapping $u \mapsto g_i(x,u)$, $i \in \{1,\ldots,p\}$, is convex on $\mathbb{R}^n$.
\end{description}
\begin{proposition}
	Let the dynamics be given by~\eqref{sys:linear} and Assumption~A5 hold. Then, the viability kernel is compact, convex, and contains the origin in its interior.
\end{proposition}
\begin{proof}
	The proof that the viability kernel is bounded, convex, and contains the origin in its interior easily adapts from the discrete-time case~\cite{BOCCIA2014}, see Propositions 5 and 6 of the appendix. As detailed in~\cite[Prop. 3.1]{Esterhuizen_Levine_arXiv_2015} and~\cite{ESTERHUIZEN_2016}, to have closedness of the viability kernel one needs to impose assumptions on the dynamics $f$, in addition to Assumption~A5. These are that $f$ is at least $C^2$ from $\mathbb{R}^n \times U_1$ to $\mathbb{R}^n$, $U_1$ an open subset containing $U$; that there exists a constant $0 < C < \infty$ such that $\sup_{u\in U} |x^T f(x,u)| \leq C (1 + \Vert x \Vert ^2)$, for all $x\in\mathbb{R}^n$; and that the set $f(x,U)$ is convex for all $x\in\mathbb{R}^n$. We note that all three of these assumptions are satisfied by the system \eqref{sys:linear}.
\end{proof}
\begin{remark}
	If the functions $g_i$ do not depend on the input, that is $g_i:\mathbb{R}^n \rightarrow \mathbb{R}$, %
	the set~$\A$ is closed without assuming convexity of $g_i(x,\cdot)$, see~\cite[Prop. 4.1]{DeDona_siam} for details.
\end{remark}
\begin{remark}
	Note that in the discrete-time setting with dynamics $x_{k+1} = f(x_k,u_k)$, a set~$S$ is said to be control invariant %
	if, for each~$x \in S$, there exists $u \in U(x)$ such that $f(x,u) \in S$ holds. Hence, establishing closedness of~$\A$ requires fewer assumptions: %
	Briefly, if $g$ is continuous, closedness can be directly shown by considering the limit of converging sequences $\{x_k\}_{k\in\mathbb{N}} \subseteq \A$.
\end{remark}

Assumption~A4 implies existence of a matrix $F$ such that the state feedback $\mu_F(x) = Fx$ renders the matrix $(A + BF)$ Hurwitz. Thus, see, e.g.\ \cite{Brockett_1970}, there exist constants $\Gamma>0$ and $\eta>0$ such that for all $x_0\in\mathbb{R}^n$: 
\begin{equation}\label{eq:feedback-estimate}
	| x(t;x_0,u_F)| \leq \Gamma e^{-\eta(t-t_0)}| x_0 |\quad \forall t\in[t_0,\infty).
\end{equation} 
These facts are used to arrive at the following proposition,  which uses ideas originally proposed in \cite{GONDHALEKAR2009}.
\begin{proposition}\label{prop:LQR}	
	Consider the system \eqref{sys:linear} under Assumptions~A4 and~A5. For all $x\in\lambda \A$, with $\lambda\in[0,1)$ there exists a constant $M(\lambda)$ such that $V_{\infty}(x)\leq M(\lambda)$.
\end{proposition}
\begin{proof}
	Since the proof essentially uses the same ideas as its discrete-time analogon, we only provide a sketch and refer to Proposition~7 of the appendix for the details. The idea of the proof is to consider a convex combination of two control functions: one that keeps the state-control pair admissible for all time, and the other, given by the linear-quadratic-regulator (LQR) feedback, driving the state to the origin asymptotically. Then, upon reaching a neighbourhood of the origin (which is guaranteed due to the choice of the convex combination) switching to the LQR feedback, which drives the state to the origin while satisfying the constraints.
\end{proof}

Proposition~\ref{prop:LQR} states that, for linear systems satisfying Assumptions~A4 and~A5, %
the infinite-horizon value function is uniformly bounded on $\lambda \A$, $\lambda\in[0,1)$. %
On the one hand, for given~$\delta$, this ensures asymptotic stability of the origin w.r.t.\ the MPC closed loop with the interior of $\A$ contained in the basin of attraction (for a sufficiently large prediction horizon~$T = N\delta$). On the other hand, the set $\mathcal{O}$ (if present, see Example~\ref{ExampleSetO}) is restricted to the boundary of $\A$.

Next, we state a theorem which combines Corollary~\ref{cor:compact-subsets} with Proposition~\ref{prop:LQR} to provide a result for the LQ case.
\begin{theorem}\label{thm:LQ-case}
	Consider the system \eqref{sys:linear} and suppose that Assumptions~A4 and~A5 hold with the symmetric, quadratic, and positive definite stage cost
		\begin{align}\label{eq:QuadraticStageCost}
			\ell(x,u) &:= (x^T, u^T)
				\left( 
					\begin{array}{cc}
						Q & N \\
						N^T & R
					\end{array}
				\right)
				\left( 
					\begin{array}{c}
						x\\
						u
					\end{array}
				\right).
		\end{align}
	Let $K\subset \Int(\A)$ be compact and $\delta > 0$ be given. %
	Then, 
	the origin is asymptotically stable w.r.t. the MPC closed loop with basin of attraction~$\mathcal{S}$ containing $K$ %
	for each prediction horizon $T = N\delta$ such that $N$ satisfies Condition~\eqref{ineq-2}.
\end{theorem}
\begin{proof}
	Clearly, there exists $\lambda \in [0,1)$ such that $K \subset \lambda \A$. Hence, Proposition~\ref{prop:LQR} implies boundedness of $V_T$ on $K$. %
	Consequently, $K\subseteq V_{\infty}^{-1}[0,\infty)\setminus \mathcal{O}$ holds. %
	Hence, Theorem~\ref{thm:main} and Corollary~\ref{cor:compact-subsets} imply the assertion supposing that Assumptions A1-A3 hold. 

	Validity of Assumptions~A1 and~A2 can be shown analogously to the discrete-time case, cp.\ \cite[Thm. 13]{BOCCIA2014}: %
	Since the cost is quadratic and positive definite, Assumption~A1 holds. In particular, we have $\ell^\star(x) = x^T Q x$,  and $\sigma_{\min}(Q)|x|^2 \leq\ell^\star(x)\leq \sigma_{\max}(Q) |x|^2$ where $\sigma_{\min}(Q)$ and $\sigma_{\max}(Q)$ denote the minimal and maximal eigenvalue of the matrix~$Q$, respectively. %
	Invoking Assumptions~A4 and~A5 imply that the unique, positive definite solution~$P$ of the algebraic Riccati equation satisfies $V_{\infty}(x) = x^T P x$ %
	on a neighbourhood~$\mathcal{N}$ of the origin. Hence, we have
	\begin{equation}
		V_\infty(x) \leq \sigma_{\max}(P) |x|^2 \leq \frac {\sigma_{\max}(P)}{\sigma_{\min}(Q)} \ell^{\star}(x) \nonumber
	\end{equation}		
	for all $x\in\mathcal{N}$, i.e.\ Assumption~A2.
	
	Next, we prove Assumption~A3, which trivially holds in the discrete-time setting with $\bar{C} = 1$. %
	To this end, we extend the previously presented argumentation based on the algebraic Riccati equation. %
	Taking the constraints into account yields $V_\delta(x) \geq x^\top  P x \geq \sigma_{\min} (P) |x|^2$. %
	Hence, Assumption~A3 holds with $\bar{C} := T \sigma_{\max}(Q) / \sigma_{\min}(P)$. 
\end{proof}
 
\begin{remark}
	The proof of Theorem~\ref{thm:LQ-case} identifies some of the constants appearing in Theorem~\ref{thm:main} for the LQ case: $\gamma = \frac{\sigma_{\max}(P)}{\sigma_{\min}(Q)}$, $\bar{C} = \delta\frac{\sigma_{\max}(Q)}{\sigma_{\min}(P)}\leq T\frac{\sigma_{\max}(Q)}{\sigma_{\min}(P)}$.
\end{remark}
	

Finally, we state the continuous-time analogon of \cite[Cor.~15]{BOCCIA2014}, which relates the sufficient horizon length to the distance of the state to $\A$'s boundary, see the appendix for a detailed proof.
\begin{corollary}\label{cor:horizon-length}
	Consider the system \eqref{sys:linear} under Assumptions~A4 and~A5, with the quadratic stage cost~\eqref{eq:QuadraticStageCost}, and a compact set $K\subset \Int(\A)$. The upper bound of the infinite-horizon value function on~$K$ is inversely proportional to the distance of the state from the boundary of the viability kernel. That is, there exists a constant $D$, such that:
	\[
		\sup_{x\in K} V_{\infty}(x) \leq \frac{D}{\dist(K;\partial \A)}.
	\]
	Moreover, for $\beta$ set to $\max\{\frac{C}{M\dist(K;\partial \A)}, \gamma\}$, if the time shift $\delta$ is smaller than $\min \{ C/M, \gamma\sqrt{\bar{C}} \}$ and $\bar{N}(K,\delta)$ satisfies~\eqref{eq:ineq_2}, the origin is asymptotically stable w.r.t. the MPC closed loop with basin of attraction~$\mathcal{S}$ containing the set~$K$.
\end{corollary}

In the future we intend to further explore conditions under which the horizon length does not blow up as the state approaches the boundary of $\A$ by using the theory of barriers as mentioned {in Section~\ref{sec:setting}}. 


\section{Example}\label{sec:example}

\begin{figure}[htb]
	\begin{center}\
		\includegraphics[width=0.7\columnwidth]{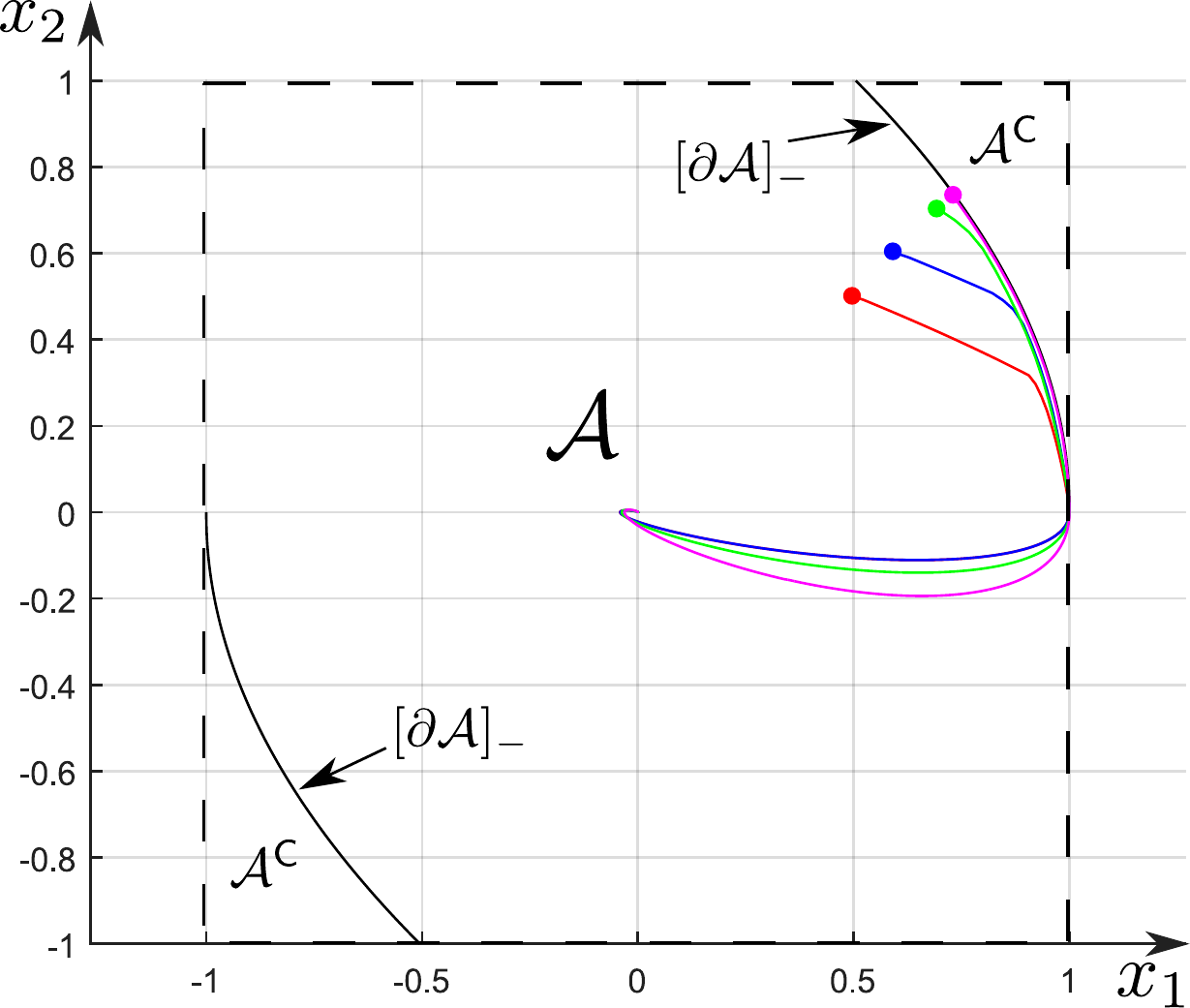} 
		\caption{Simulation of the MPC closed loop from various initial conditions with $\delta = 0.1$. Red curve: $x_0 = (0.5, 0.5)^T$, $N = 4$; blue curve: $x_0 = (0.6, 0.6)^T$, $N = 4$; green curve: $x_0 = (0.7, 0.7)^T$, $N = 5$; magenta curve: $x_0 = (0.733, 0.73)^T$, $N = 7$.} 
		\label{fig:double_integrator}                             
	\end{center}                                
\end{figure}

\noindent We demonstrate the growth of the sufficient horizon length $N$ for a chosen time shift~$\delta$ as the initial condition approaches the viability kernel's boundary. We consider the double integrator:
\begin{align*}
	\dot{x}_1(t) & = x_2(t), \qquad \dot{x}_2(t) = u(t),
\end{align*}
with $|u| \leq 1$ and $|x_i| \leq 1$, {$i \in \{1,2\}$}. %
As shown in \cite{DeDona_siam}, parts of the boundary of the viability kernel (called the barrier) consist of integral curves that, together with a particular control function, satisfy a minimum-like principle and %
intersect the boundary of the constrained state space tangentially. %
We use this fact to construct the two solid curves that form the barrier, labelled $\DAM$, see Figure~\ref{fig:double_integrator}. %

We run the continuous-time MPC algorithm (solving the optimisation problem via a direct method) with the running cost $\ell(x,u) = x_1^2 + x_2^2 + u^2$ from various initial conditions that approach the boundary of $\A$, with various time shifts~$\delta$ and horizons~$N$. Table~\ref{tab:table1} displays the smallest $N = N(x_0,\delta)$, for which MPC steers the particular initial state to the origin while maintaining constraint satisfaction. It is interesting to note that the magenta curve, initiating from $x_0\in \DAM$, results from a finite horizon length, $N = 7$.  This emphasises that the infinite-horizon value function may be bounded on the viability kernel's boundary, and that the statement in Corollary~\ref{cor:horizon-length} is merely sufficient.

\begin{table}[h!]
	\begin{center}
		\caption{Smallest horizon length $N(\delta)$ for time shift $\delta$, from $x_0$.}
		\begin{tabular}{l|c|c|c}
			\toprule 
			$x_0$ & $\delta = 0.1$ & $\delta = 0.05$ & $\delta = 0.03$ \\		
			\midrule 
			$(0.5,0.5)^T$ & 4 & \phantom{0}7 & 10\\ \midrule %
			$(0.6,0.6)^T$ & 4 & \phantom{0}7 & 11\\ \midrule %
			$(0.7,0.7)^T$ & 5 & 10 & 14\\ \bottomrule %
		\end{tabular}
		\label{tab:table1}
	\end{center}
\end{table}

%
%
%
%
%
%
%
%
%

\vspace{-0.5cm}
\section{Conclusion}\label{sec:conclusion}

\noindent {We considered MPC without (stabilising) terminal costs and/or constraints for continuous-time systems. %
We proposed sufficient conditions to ensure asymptotic stability of an 
equilibrium. 
Simultaneously, we showed recursive feasibility 
for all 
initial values contained in some compact set. %
In particular, we provided a formula that allows one to determine a prediction horizon such that the finite-horizon value function is guaranteed to enjoy a Lyapunov-like decrease along trajectories of the MPC closed loop. %
For the linear-quadratic case, the interior of the viability kernel can (essentially) be covered using this technique using standard assumptions.}

\section*{ACKNOWLEDGEMENT}

\noindent This work was 
funded by BMBF-Projekts 05M18OCA / 05M18SIA within KONSENS: ``Konsistente Optimierung und Stabilisierung elektrischer Netzwerksysteme''. Karl Worthmann thanks the German Research Foundation (DFG) for their support (Grant No.\ WO\,2056/6-1).

\bibliographystyle{IEEEtran}
\bibliography{bibliography}

\begin{thebibliography}{10}
\providecommand{\url}[1]{#1}
\csname url@samestyle\endcsname
\providecommand{\newblock}{\relax}
\providecommand{\bibinfo}[2]{#2}
\providecommand{\BIBentrySTDinterwordspacing}{\spaceskip=0pt\relax}
\providecommand{\BIBentryALTinterwordstretchfactor}{4}
\providecommand{\BIBentryALTinterwordspacing}{\spaceskip=\fontdimen2\font plus
\BIBentryALTinterwordstretchfactor\fontdimen3\font minus
  \fontdimen4\font\relax}
\providecommand{\BIBforeignlanguage}[2]{{%
\expandafter\ifx\csname l@#1\endcsname\relax
\typeout{** WARNING: IEEEtran.bst: No hyphenation pattern has been}%
\typeout{** loaded for the language `#1'. Using the pattern for}%
\typeout{** the default language instead.}%
\else
\language=\csname l@#1\endcsname
\fi
#2}}
\providecommand{\BIBdecl}{\relax}
\BIBdecl

\bibitem{RawlMayn18}
{J. B. Rawlings, D. Q. Mayne, M. M. Diehl}, \emph{Model Predictive Control:
  Theory, Computation, and Design}, 2nd~ed.\hskip 1em plus 0.5em minus
  0.4em\relax Nob Hill Publishing, 2018.

\bibitem{NeviPrim97}
V.~Nevisti{\'c} and J.~A. Primbs, ``Receding horizon quadratic optimal control:
  performance bounds for a finite horizon strategy,'' in \emph{Proc. European
  Control Conf. (ECC)}, 1997, pp. 3584--3589.

\bibitem{JadbYu01}
A.~Jadbabaie, J.~Yu, and J.~Hauser, ``Unconstrained receding-horizon control of
  nonlinear systems,'' \emph{IEEE Trans. Automat. Control}, vol.~46, no.~5, pp.
  776--783, 2001.

\bibitem{GrimMess05}
{G. Grimm, M. J. Messina, S. E. Tuna, A. R. Teel}, ``Model predictive control:
  for want of a local control {L}yapunov function, all is not lost,''
  \emph{IEEE Trans. Automat. Control}, vol.~50, no.~5, pp. 546--558, 2005.

\bibitem{GrunPann10}
L.~Gr{\"u}ne, J.~Pannek, M.~Seehafer, and K.~Worthmann, ``Analysis of
  unconstrained nonlinear {MPC} schemes with varying control horizon,''
  \emph{SIAM J. Control Optim.}, vol.~48, no.~8, pp. 4938--4962, 2010.

\bibitem{REBLE2012}
M.~Reble and F.~Allg{\"o}wer, ``Unconstrained model predictive control and
  suboptimality estimates for nonlinear continuous-time systems,''
  \emph{Automatica}, vol.~48, no.~8, pp. 1812 -- 1817, 2012.

\bibitem{Worthmann2014}
K.~Worthmann, M.~Reble, L.~Gr\"{u}ne, and F.~Allg\"{o}wer, ``The role of
  sampling for stability and performance in unconstrained nonlinear model
  predictive control,'' \emph{SIAM J. Control Optim.}, vol.~52, no.~1, pp.
  581--605, 2014.

\bibitem{cannon_2015}
M.~S. Darup and M.~Cannon, ``{A missing link between nonlinear MPC schemes with
  guaranteed stability},'' in \emph{Proc. 54th IEEE Conf. on Decision and
  Control (CDC)}, 2015, pp. 4977--4983.

\bibitem{GrunPann17}
L.~Gr\"{u}ne and J.~Pannek, \emph{Nonlinear Model Predictive Control. Theory
  and Algorithms}, 2nd~ed.\hskip 1em plus 0.5em minus 0.4em\relax London:
  Springer, 2017.

\bibitem{BOCCIA2014}
A.~Boccia, L.~Gr{\"u}ne, and K.~Worthmann, ``Stability and feasibility of state
  constrained {MPC} without stabilizing terminal constraints,'' \emph{Systems
  \& Control Letters}, vol.~72, pp. 14 -- 21, 2014.

\bibitem{boccia_2014_linear}
------, ``{Stability and feasibility of state-constrained linear MPC without
  stabilizing terminal constraints },'' in \emph{Proc. 21st Int. Symp. Math.
  Theory Networks Syst.}, 2014, pp. 453--460.

\bibitem{FaulBonv15b}
T.~Faulwasser and D.~Bonvin, ``On the design of economic {NMPC} based on an
  exact turnpike property,'' \emph{IFAC-PapersOnLine}, vol.~48, no.~8, pp.
  525--530, 2015.

\bibitem{FaulBonv15}
------, ``On the design of economic {NMPC} based on approximate turnpike
  properties,'' in \emph{Proc. 54th IEEE Conf. on Decision and Control (CDC)},
  2015, pp. 4964--4970.

\bibitem{Coron_2019}
J.-M. Coron, L.~Gr{\"u}ne, and K.~Worthmann, ``Model predictive control, cost
  controllability, and homogeneity,'' \emph{arXiv:1906.05112}, 2019.

\bibitem{DeDona_siam}
J.~De~Dona and J.~L{\'e}vine, ``On barriers in state and input constrained
  nonlinear systems,'' \emph{SIAM J. Control Optim.}, vol.~51, no.~4, pp.
  3208--3234, 2013.

\bibitem{aubin2009viability}
J.-P. Aubin, \emph{Viability theory}.\hskip 1em plus 0.5em minus 0.4em\relax
  Springer Science \& Business Media, 2009.

\bibitem{sontag2013}
E.~D. Sontag, \emph{Mathematical control theory: deterministic finite
  dimensional systems}.\hskip 1em plus 0.5em minus 0.4em\relax Springer Science
  \& Business Media, 2013, vol.~6.

\bibitem{DaruWort18}
M.~S. Darup and K.~Worthmann, ``{Tailored MPC for mobile robots with very short
  prediction horizons},'' in \emph{Proc. European Control Conf. (ECC)}, 2018,
  pp. 1361--1366.

\bibitem{Wort11}
\BIBentryALTinterwordspacing
K.~Worthmann, ``Stability analysis of unconstrained receding horizon control
  schemes,'' Ph.D. dissertation, University of Bayreuth, 2011. [Online].
  Available: \url{https://epub.uni-bayreuth.de/273/}
\BIBentrySTDinterwordspacing

\bibitem{WortMehr15}
K.~Worthmann, M.~W. Mehrez, M.~Zanon, G.~K.~I. Mann, R.~G. Gosine, and
  M.~Diehl, ``Regulation of differential drive robots using continuous time mpc
  without stabilizing constraints or costs,'' \emph{IFAC-PapersOnLine},
  vol.~48, no.~23, pp. 129--135, 2015.

\bibitem{saintpierre_1994}
P.~Saint-Pierre, ``Approximation of the viability kernel,'' \emph{Applied
  Mathematics and Optimization}, vol.~29, no.~2, pp. 187--209, 1994.

\bibitem{monnet_2016}
D.~Monnet, J.~Ninin, and L.~Jaulin, ``Computing an inner and an outer
  approximation of the viability kernel,'' \emph{Reliable Computing}, vol.~22,
  2016.

\bibitem{Esterhuizen_Levine_arXiv_2015}
W.~Esterhuizen and J.~L{\'e}vine, ``Barriers in nonlinear control systems with
  mixed constraints,'' \emph{arXiv:1508.01708}, 2015.

\bibitem{ESTERHUIZEN_2016}
------, ``Barriers and potentially safe sets in hybrid systems: Pendulum with
  non-rigid cable,'' \emph{Automatica}, vol.~73, pp. 248 -- 255, 2016.

\bibitem{Brockett_1970}
R.~W. Brockett, \emph{Finite Dimensional Linear Systems}.\hskip 1em plus 0.5em
  minus 0.4em\relax John Wiley \& Sons, 1970.

\bibitem{GONDHALEKAR2009}
R.~Gondhalekar, J.~Imura, and K.~Kashima, ``Controlled invariant feasibility
  — a general approach to enforcing strong feasibility in {MPC} applied to
  move-blocking,'' \emph{Automatica}, vol.~45, no.~12, pp. 2869 -- 2875, 2009.

\end{thebibliography}

\section{APPENDIX}

\subsection{Some facts of the admissible set}\label{SubsectionAdmissibleSet}

We summarise some facts concerning the admissible set for a linear time-invariant system with assumptions A4 and A6. Most of these facts easily adapt from the discrete-time case, as in \cite{BOCCIA2014}.

\begin{definition}
	A set $S$ is said to be \emph{control invariant} with respect to the system \eqref{eq:system-1}, provided that for all $x_0\in S$ $\U_{\infty}(x_0)$ is nonempty and there exists a $u\in\U_{\infty}(x_0)$ such that $x(t;x_0,u)\in S$ for all $t\in[t_0,\infty[$.
\end{definition}
\begin{definition}
	The maximal control invariant set contained in $X$, which we label $\C$, is the union of all control invariant sets that are subsets of $X$.
\end{definition}

\begin{proposition}\label{prop:A-contr-invariant}
	$\A = \C$.
\end{proposition}
\begin{proof}
	By definition, $\C\subseteq \A$. Let us show that $\A \subseteq \C$ by contradiction. Consider a point $x_0\in\A$ at $t_0$ and suppose that $\A$ is not control invariant. Then, for all $u\in\U_{\infty}(x_0)$ there exists a $t_u\in[t_0,\infty[$ such that $x_u := x(t_u;x_0,u)\notin\A$, and thus $\mathcal{U}_{\infty}(x_u) = \emptyset$. Thus $\U_{\infty}(x_0)=\emptyset$, which contradicts that fact that $x_0\in \A$. Therefore, $\A\subseteq\C$, and thus $\C = \A$.
\end{proof}

\begin{proposition}\label{PropositionConvexity}
	Consider system \eqref{sys:linear}. If $\E$ is convex, then $\A$ is convex.
\end{proposition}
\begin{proof}
	Consider any initial condition in the convex hull of the two points $x_1,x_2\in \A$, that is, consider $x_3 = \rho x_1 + (1-\rho) x_2$, $\rho\in[0,1]$. Consider the same convex combination of two admissible inputs associated with $x_1$ and $x_2$ at every $t\in[0,\infty[$, that is, consider $u_3(t) = \rho u_1(t) + (1-\rho) u_2(t)$, with $u_i\in\mathcal{U}_{\infty}(x_i)$, $i=1,2$. By linearity it can be verified that $x(t;x_3,u_3) = \rho x(t;x_1,u_1) + (1 - \rho) x(t;x_2,u_2)$. We have $(x(t;u_i,x_i),u_i(t))\in \E$ for all $t\in[t_0,\infty[$, $i=1,2$. Because $\E$ is convex, we have $(x(t;x_3,u_3),u_3(t)) \in \E$ for all $t\in[t_0,\infty[$, and thus $x_3\in\A$.
\end{proof}
Clearly if $\E$ is bounded then $\A$ is bounded.
\begin{proposition}\label{prop:0-in-A}
	For the system \eqref{sys:linear} with assumptions A4 and A6 the origin is in the interior of $\A$.
\end{proposition}
\begin{proof}
	Clearly $0\in\A$. We need to show that there exists a neighbourhood of $0$ contained in $\A$. From \eqref{eq:feedback-estimate}, with a stabilizable pair $(A,B)$, we have $| x(t;x_0,u_F)| \leq \Gamma e^{-\eta(t-t_0)}| x_0 |$ for all $t\geq t_0$. Moreover, $| Fx(t;x_0,u_F) | \leq \Vert F \Vert \Gamma e^{-\eta(t-t_0)}| x_0 |$ for all $t\geq t_0$. Thus, for any $\ee > 0$ there exists a $\nu := \frac{\ee}{(1 + \Vert F \Vert)\Gamma}$ such that for all $(x_0, Fx_0)\in \nu\mathcal{B}_{n+m}$ we have $(x(t;x_0,u_F),Fx(t;x_0,u_F))\in \ee \mathcal{B}_{n+m}$, where $\mathcal{B}_{n+m}\subset\mathbb{R}^{n+m}$ is the open unit ball about the origin of $\E$. Recall that $0\in \Int(\E)$. Thus, we can select $\ee>0$ small enough such that $\ee \mathcal{B}_{n+m}\subset \E$, and so $\mathsf{proj}_{\mathbb{R}^n}(\nu \mathcal{B}_{n+m})\subset \A$.
\end{proof}

\begin{proposition}\label{prop:lambdaA-invar}
	For the system \eqref{sys:linear} with assumptions A4 and A6 the set $\lambda\A$ is control invariant for any $\lambda\in[0,1]$.
\end{proposition}
\begin{proof}
	Recall that $0\in \A$ (Proposition~\ref{prop:0-in-A}) and that $\A$ is control invariant (Proposition~\ref{prop:A-contr-invariant}). Consider a point $x_0\in\lambda \A$, then $x_0/ \lambda \in \A$. Thus, there exists a $u\in\U_{\infty}(x_0/\lambda)$ such that $x(t;x_0/ \lambda,u)\in \A$ for all $t\geq t_0$. By linearity it may be verified that $\lambda x(t;x_0/\lambda, u) = x(t;x_0,\lambda u)\in\lambda \A$ for all $t\geq t_0$, with $\lambda u\in\mathcal{U}_{\infty}(x_0)$ (from A6). Thus, $\lambda u$ renders $\lambda\A$ invariant.
\end{proof}

\subsection{Proof of Proposition~\ref{prop:LQR}}\label{sec:app-proof-prop2}

Consider $x_0\in\lambda\A$, for which there exists a control, which we label $u_{\lambda}\in\mathcal{U}_{\infty}(x_0)$, such that $x(t;x_0, u_{\lambda})\in\lambda \A $ for all $t\geq t_0$ (Proposition~\ref{prop:lambdaA-invar}). From \eqref{eq:feedback-estimate}, because $(A,B)$ is stabilizable, there exists a $u_{F}\in\U_{\infty}(x_0)$ such that  $| x(t;x_0,u_F)|\leq \Gamma e^{-\eta(t-t_0)}| x_0 |$. Thus, $| (x(t;x_0,u_F), u_F(t))|\leq L\lambda d_{\min}$, where $L: = (1 + \Vert F \Vert)\Gamma d_{\max}d^{-1}_{\min}$, $d_{\max}:=\sup_{x\in X} | x |$, $d_{\min}:=\inf_{x\in\partial X} | x |$. If $\lambda L \leq 1$ then the solution remains in $\lambda L\E$ for all time and converges to the origin, and $V_{\infty}(x_0) = \sup_{x\in\lambda\A} J_{\infty}(x,u_F)\leq \alpha$, $\alpha \geq 0$. Otherwise, if $\lambda L>1$, consider $\mu x(t;x_0, u_{\lambda}) + (1 - \mu) x(t;x_0,u_F)$, with $\mu\in[0,1]$, which is the solution obtained with the control $\tilde{u} := \mu u_{\lambda} + (1 - \mu) u_F$. It may be verified that: $| (x(t;x_0,\tilde{u}), \tilde{u}(t))| \leq  [\mu\lambda + (1 - \mu) L \lambda] d_{\min}$. Choose $\mu$ such that $\mu\lambda + (1 - \mu)L \lambda = 1$, then $(x(t;x_0,\tilde{u}), \tilde{u}(t))\in \E$ for all time. Note that $\mu > 0$, (if $\mu = 0$ we would have $u = u_F$ and $\lambda L \leq 1$.) With this control, there exists a $\bar{t}\geq t_0$ such that $x(\bar{t};x_0,u)\in\ee\lambda\A$, where we have defined $\ee:= 1 - \frac{1 - \lambda}{\lambda L}$. From our choice of $\mu$ we also see that $\ee\in]\mu,1[$. If we consider time $m\bar{t}$, with $m\in\mathbb{N}_0$, we get $x(m\bar{t};x_0,u)\in \ee^m \lambda\A$. Let $m(x_0)$ be the smallest integer such that $\ee^{m(x_0)} \lambda L <1$, and let $\tilde{x} := x(m(x_0)\bar{t};x_0,u)$. If we now switch to the feedback $u_F$, we get  $| x(t;\tilde{x},u_F)| \leq \ee^{m(x_0)} \lambda L d_{\min}$ for all $t\geq m(x_0)\bar{t}$. Thus, $x(t;\tilde{x},u_F)\in \lambda L\ee^{m(x_0)} X$ for all $t\geq m\bar{t}$, and converges to the origin. We get that:
\[
V_{\infty}(x_0) = J_{m(x_0)\bar{t}}(x_0,u) + J_{\infty}(\tilde{x},u_F) \leq m(x_0)\bar{t}\beta + \alpha < \infty.
\]

\subsection{Proof of Corollary~\ref{cor:horizon-length}}

The proof adapts from \cite[Cor. 14]{BOCCIA2014}. As in the proof of Proposition~\ref{prop:LQR}, consider the constants $\ee$ and $L$ along with the mapping: $m(x) := \inf\{m\in\mathbb{N}: \ee^{m}\leq \frac{1}{\lambda L}\}$. Then it can be shown, after conducting an asymptotic analysis, that $m(x) \thicksim \frac{wL\ln(L)}{\dist(x;\partial \A)}$ as $x\rightarrow \partial\A$, where $w \in[\inf_{x\in\partial\A} |x|, \sup_{x\in\A} |x|]$. The result then follows from the fact, as established in the proof of Proposition~\ref{prop:LQR}, that
$V_{\infty}(x)  \leq m(x)\bar{t}\beta + \alpha$.

\end{document}